\documentclass{amsproc}

\usepackage{enumitem}
\usepackage{tikz}
\usepackage{tikz-cd}
\usepackage{amssymb}
\usepackage{blkarray}
\usepackage{mathdots}
\usepackage{hyperref}

\DeclareMathOperator{\ad}{\text{ad}}
\DeclareMathOperator{\Tr}{\text{Tr}}
\DeclareMathOperator{\id}{\text{id}}
\DeclareMathOperator{\Span}{\text{Span}}
\DeclareMathOperator{\Sym}{\text{Sym}}
\newcommand{\M}{\mathbb{M}}
\newcommand{\projh}[1]{\overline{#1}} 
\newcommand{\projq}[1]{\widehat{#1}} 
\DeclareMathOperator{\m}{\ell} 
\newcommand{\heis}{\mathcal{H}}
\newcommand{\V}{\mathcal{V}}
\newcommand{\setdif}{-}
\newtheorem{theorem}{Theorem}[section]
\newtheorem{lemma}[theorem]{Lemma}
\newtheorem{corollary}[theorem]{Corollary}
\newtheorem{proposition}[theorem]{Proposition}

\theoremstyle{definition}

\theoremstyle{remark}

\usetikzlibrary{decorations.markings,arrows}

\pgfdeclarelayer{edgelayer}
\pgfdeclarelayer{nodelayer}
\pgfsetlayers{edgelayer,nodelayer,main}

\tikzstyle{none}=[]
\tikzstyle{root1}=[fill={rgb,255: red,120; green,120; blue,120}, draw=black, shape=circle]
\tikzstyle{root2}=[fill={rgb,255: red,45; green,45; blue,45}, draw=black, shape=circle]
\tikzstyle{root3}=[fill={rgb,255: red,220; green,220; blue,220}, draw=black, shape=circle]
\tikzstyle{horizontal_tick}=[fill=white, draw=black, shape=rectangle]

\tikzstyle{plainedge}=[-, draw=black=]

\numberwithin{equation}{section}

\begin{document}

\title{The structure of non-Fricke Monstrous Lie algebras}

\author{Daniel Tan}
\address{Department of Mathematics, Rutgers University, 110 Frelinghuysen Rd, Piscataway, New Jersey 08854}
\curraddr{}
\email{dwt24@math.rutgers.edu}
\thanks{I would like to thank Professor Lisa Carbone, Hong Chen, Jishen Du, Dennis Hou and Forrest Thurman for working on this project. I would also like to thank the AMS-UMI Joint Meeting in Palermo for the opportunity to speak.}


\date{}

\begin{abstract}
We consider Borcherds algebras with no real roots and the property that all zeroes in the Borcherds Cartan matrix occur in a single diagonal zero block. It follows that all other entries of the matrix are negative. We give a structure
theorem for these Borcherds algebras, decomposing them into
free, Heisenberg and abelian subalgebras. We show that a class of such Borcherds algebras are the Monstrous Lie algebras associated to non-Fricke elements of the Monster finite simple
group. This new perspective on their structure gives an efficient method to
compute their twisted denominator formulas.

\end{abstract}

\maketitle

\section{Introduction}
Let $ \M $ be the Monster finite simple group. Let $ V^\natural $ be the Moonshine module vertex operator algebra constructed and shown to have automorphism group $\mathbb{M}$ by Frenkel, Lepowsky and Meurman \cite{FLM}. The Conway--Norton Monstrous Moonshine Conjecture \cite{CN} was partially settled by Frenkel, Lepowsky and Meurman using the Moonshine module $V^\natural $ in a central role. The  Monster Lie algebra $\mathfrak m$, constructed by Borcherds \cite{BoInvent},  was a crucial ingredient in the proof of the remaining cases of Monstrous Moonshine.  

Borcherds discovered $\mathfrak{m}$ as an example of a  new class of Lie algebras, called Borcherds  algebras or generalized Kac--Moody algebras  \cite{BoGKMA}. He defined $\mathfrak{m}$ as a certain subquotient of the vertex algebra $V=V^\natural\otimes V_{II_{1,1}}$, where $V_{II_{1,1}}$ is the vertex algebra associated to the 2-dimensional unimodular even  Lorentzian lattice  $II_{1,1}$.

Using the No-ghost theorem \cite{BoInvent}, Borcherds gave $\mathfrak{m} $ a $ II_{1,1} $-grading into root spaces and determined the dimensions of the root spaces in terms of the dimensions of the weight-homogeneous subspaces of $ V^\natural $, that is, the coefficients of the normalized modular $ j $-function $j(\tau) - 744 $. With the dimensions of the root spaces, he showed that $ \mathfrak{m} $ has a realization $ \mathfrak{m}=\mathfrak{g}(A)/\mathfrak{z} $ in terms of generators and relations, where $A$ is a Borcherds Cartan matrix related to $ j(\tau) - 744 $ and $\mathfrak{z}$ is the center of the Borcherds algebra $\mathfrak{g}(A)$. The No-ghost theorem also transfers the action of $\M$ on $ V^\natural $ to a $ II_{1,1} $-grading preserving action of $\M$ on $ \mathfrak{m} $. Borcherds derived `twisted' versions of the denominator formula for $ \mathfrak{m} $ by taking the trace of the action of each $ g \in \M $ in the usual denominator formula for $ \mathfrak{m}$. These twisted denominator formulas imply certain relations between the coefficients of the McKay--Thompson series $ T_g(q) $. Borcherds then extended the homological methods of Garland and Lepowsky  \cite{GarLep} to derive the replication formulas conjectured by Conway and Norton \cite{CN}.

In \cite{JurJPAA}, Jurisich obtained a decomposition of the Monster Lie algebra  into two free subalgebras $ \mathfrak{u}^\pm $ and a copy of $ \mathfrak{gl}_2 \cong \mathfrak{sl}_2 \oplus \mathbb{C}h $. Furthermore, the two free subalgebras are freely generated by $ \mathfrak{sl}_2 $-modules, given by the adjoint action. The proof of the existence of such a decomposition uses the explicit form of the Borcherds Cartan matrix for $\mathfrak m$, also given in \cite{JurJPAA}. The homology of the positive free part $ \mathfrak{u}^+ $ together with $ \mathfrak{sl}_2 $-representation theory was used in \cite{JLW} to efficiently compute the twisted denominator formulas for the Monster Lie algebra.

In \cite{No}, Norton  generalized the Conway--Norton  Monstrous Moonshine Conjecture. The generalized conjecture involves the unique irreducible $g$-twisted $V^\natural$-module $V^\natural_g$, shown to exist by Dong, Li and Mason \cite{DLiM}, for each $ g \in \M $. Furthermore, associated to each conjugacy class $ [g] $ of $ \M $ with representative $ g $, is an abelian intertwining algebra structure on the direct sum of $V^\natural_{g^i}$ for $ i = 0, \dots, o(g)-1 $, which was shown to exist in \cite{vEMS}, but believed to exist earlier. They also showed that the abelian intertwining algebra associated to $ [g] $ has a projective action of the centralizer $ C_{\M}(g) $. 

 As part of his work on the Generalized Moonshine Conjecture \cite{CarThesis, CarGMI, CarDuke}, Carnahan  constructed a family of Lie algebras $\mathfrak m_g$, one for each conjugacy class $ [g] $ of $\mathbb{M}$. He constructed $\mathfrak m_g$ as a certain subquotient of the abelian intertwining algebra $ \bigoplus_{i=0}^{o(g)-1} V^\natural_{g^i} $ tensored with $ V_{II_{1,1}} $. Using the No-ghost theorem, he showed that $ \mathfrak{m}_g $ has a Borcherds Cartan matrix description \cite{CarDuke} with simple root space dimensions related to the McKay--Thompson series $ T_g(\tau) $. Furthermore, he showed that the projective action of $ C_{\M}(g) $ on the abelian intertwining algebra is transferred to a root space preserving action on $ \mathfrak{m}_g $. Carnahan \cite{CarThesis} computed the twisted denominator formulas for $ \mathfrak{m}_g $ twisted by elements in a central extension of $ C_{\M}(g) $ using homological methods similar to those in \cite{GarLep}.

We recall that an element $g\in\mathbb{M}$ is called {\it Fricke} if  the McKay--Thompson series $T_{g}(\tau)$ is invariant
under the level $N$ Fricke involution $\tau \mapsto -1/N\tau$ for some $N\geq 1$. Otherwise $g$ is called non-Fricke.

In the case where $g$ is Fricke, Carnahan \cite{CarFricke} proved that $ \mathfrak{m}_g $ decomposes analogously to Jurisich's decomposition $ \mathfrak m \cong\mathfrak u^-\oplus( \mathfrak{sl}_2 \oplus \mathbb{C}h )\oplus\mathfrak u^+$, where $\mathfrak u^\pm$ are free Lie algebras. 
The summand $ \mathfrak{sl}_2 $ occurs because the Fricke type Monstrous Lie algebras have exactly one real simple root. Carnahan also showed that the Borcherds Cartan matrix for $\mathfrak{m}_g$ has the same form as that of $\mathfrak{m}$, with different block sizes as given by the dimensions of the root spaces for imaginary simple roots of $\mathfrak{m}_g$.
With this decomposition and homology computations similar to those in \cite{JLW}, Carnahan greatly simplified the computation of the twisted denominator formulas for $\mathfrak{m}_g$. 

Our objective in this work is to prove a similar decomposition for $\mathfrak{m}_g$ when $g$ is non-Fricke. Our motivation, in part, is to simplify the computation of the twisted denominator formulas for $\mathfrak{m}_g$.  And furthermore, the existence of free subalgebras has led to other applications, such as constructions of associated Lie group analogs \cite{CJM}. For Fricke-type $g $, Carbone and Jurisich use the free subalgebras to associate a Magnus group construction to infinite-dimensional Lie algebra $ \mathfrak{m}_g $ \cite{CJ}.

To prove an analog for $\mathfrak m_g$ of Jurisich's decomposition  $\mathfrak m\cong\mathfrak u^-\oplus (\mathfrak{sl}_2 \oplus \mathbb{C}h )\oplus \mathfrak u^+$ when $g$ is non-Fricke, we  consider general Borcherds algebras with no real roots and the property that all zeroes in the Borcherds Cartan matrix occur in a single diagonal block. It follows that all other entries of the matrix are negative. We give a structure theorem for this class of Borcherds algebras in Section~\ref{section:general decomposition}, decomposing them as $\mathfrak{u}^-\oplus (\mathcal{H} \oplus \mathbb{C}h)\oplus \mathfrak{u}^+$, where $\mathfrak{u}^\pm$ are free subalgebras generated by certain imaginary root vectors and  $\mathcal H $ is an infinite-dimensional Heisenberg Lie algebra.

When $g$ is non-Fricke, we give the explict Borcherds Cartan matrices for $\mathfrak m_g$ in Section~\ref{section:monstrous Lie algebras} as sketched by Carnahan in \cite{CarDuke}. The block sizes are given by the dimensions of the simple imaginary root spaces. These in turn  are related to the McKay--Thompson series $ T_g(\tau) $. We then show that the $\mathfrak m_g$ are examples of such a class of  Borcherds algebras studied in Section~\ref{section:general decomposition}. 

The non-Fricke Monstrous Lie algebras have no real roots, hence the $\mathfrak{sl}_2$ summand that we see in the Fricke case does not occur. The analogous summand in this case is the infinite-dimensional Heisenberg Lie algebra $\mathcal H $. This decomposition,  together with $\mathcal H $-representation theory, allows us to efficiently compute the twisted denominator formulas Section~\ref{section:twisted denominator}.

This work is a report on a larger project \cite{CCDHTT} where we study the structure of $ \mathfrak{m}_g $ as Borcherds algebras when $ g \in \M $ is non-Fricke. In the larger project, we will show the existence of the sequences $ \{ c(m,0) \}_{m = 1}^\infty $ and $ \{ c(1, n/N) \}_{n=1}^\infty $ from Section \ref{section:monstrous Lie algebras} and the explicit Cartan matrices will be constructed in detail.

\section{Borcherds algebras} \label{section:GKMs}

In this section, we briefly recall some key definitions and properties about Borcherds algebras, first given in \cite{BoGKMA}. The proofs of these results can be found in \cite{jurisich1996exposition}.

Let $ I $ be a countable index set and let $ A = (a_{ij})_{i,j \in I} $ be a real matrix. 
We say that $ A $ is a \emph{symmetrized Borcherds Cartan matrix} if it satisfies the following conditions:
\begin{enumerate}[label=(C\arabic*)]
\item The matrix $ A $ is symmetric.
\item If $ i \neq j $, then $ a_{ij} \leq 0 $.
\item If $ a_{ii} > 0 $, then $ 2a_{ij}/a_{ii} \in \mathbb{Z} $ for all $ j \in I $. 
\end{enumerate}

The \emph{Borcherds algebra $ \mathfrak{g}(A) = \mathfrak{g} $ associated to the matrix $ A $} is Lie algebra with generators $ \{ e_i, f_i, h_i\}_{i \in I} $ and relations given by the ideal generated by
\begin{enumerate}[label=(R\arabic*)]
\item $ [h_i,h_j] $, $ [h_i,e_j] - a_{ij} e_j $, $ [h_i, f_j] + a_{ij} f_j $, $ [e_i,f_j] - \delta_{ij} h_i $, for all $ i,j \in I $; \label{item:g_0 relations}
\item $ (\ad e_i)^{1 - 2 a_{ij}/a_{ii}} e_j $, $ (\ad f_i)^{1 - 2 a_{ij}/a_{ii}} f_j $, for all $ i,j \in I $ with $ i \neq j $, $ a_{ii} > 0 $; \label{item:g relations 1}
\item $ [e_i,e_j] $ , $ [ f_i, f_j] $, for all $ i,j \in I $ with $ a_{ij} = 0 $. \label{item:g relations 2}
\end{enumerate}
Any quotient of $ \mathfrak{g}(A) $ by an ideal contained in its center is a \emph{Borcherds algebra}. Borcherds algebras are also known as \emph{generalized Kac-Moody algebras}. 

We will also work with the Lie algebra $ \mathfrak{g}_0(A) = \mathfrak{g}_0 $ with generators $ \{ e_i, f_i, h_i\}_{i \in I} $ and relations given by \ref{item:g_0 relations} above. This algebra has the triangular decomposition $ \mathfrak{g}_0 = \mathfrak{n}_0^- \oplus \mathfrak{h}_0 \oplus \mathfrak{n}_0^+ $, where the subalgebra $ \mathfrak{h}_0 $ is abelian with a basis $ \{h_i\}_{i\in I} $, and the subalgebras $ \mathfrak{n}_0^\pm $ are freely generated by $ \{e_i\}_{i\in I} $ and $ \{f_i\}_{i\in I} $, respectively. The algebra $ \mathfrak{g} $ has the triangular decomposition $ \mathfrak{g} = \mathfrak{n}^- \oplus \mathfrak{h} \oplus \mathfrak{n}^+ $ with $ \mathfrak{h}_0 \cong \mathfrak{h} $ and $ \mathfrak{n}^\pm \cong \mathfrak{n}_0^\pm/ \mathfrak{k}_0^\pm $, 
where $ \mathfrak{k}_0^\pm $ is the ideal generated by the elements \ref{item:g relations 1} and \ref{item:g relations 2} in $ \mathfrak{n}_0^\pm $. 

To view the root system of $ \mathfrak{g} $ with linearly independent simple roots, we extend $ \mathfrak{g} $ with the space $ \mathfrak{d} $ spanned by the derivations $ D_i $, $ i \in I $, on $ \mathfrak{g} $ determined by
\begin{equation*}
D_i e_j = \delta_{ij} e_i, \ D_i f_j = - \delta_{ij} f_i, \ D_i h_j = 0, \text{ for all } j \in I .
\end{equation*}
Defining the \emph{extended Borcherds algebra} $ \mathfrak{g}^e = \mathfrak{d} \ltimes \mathfrak{g} $, gives the enlarged Cartan subalgebra of $ \mathfrak{g}^e = \mathfrak{n}^- \oplus \mathfrak{h}^e \oplus \mathfrak{n}^+ $ to be $ \mathfrak{h}^e = \mathfrak{h} \oplus \mathfrak{d} $.

For each $ i \in I $, let the \emph{simple root} $ \alpha_i \in (\mathfrak{h}^e)^* $ be defined by the conditions
\begin{equation}\label{eq:simple root def}
[h, e_i] = \alpha_i(h) e_i, \text{ for all } h \in \mathfrak{h}^e.
\end{equation}
With the extension, the simple roots are linearly independent. A \emph{root} of $ \mathfrak{g} $ is a non-zero element $ \phi \in (\mathfrak{h}^e)^* $ such that the space
\begin{equation}\label{eq:root space}
\mathfrak{g}_\phi = \{ x \in \mathfrak{g} : [h, x] = \phi(h) x, \text{ for all } h \in \mathfrak{h}^e \} 
\end{equation}
is non-zero. In this case, we call $ \mathfrak{g}_\phi $ the \emph{root space} of $ \mathfrak{g} $ corresponding to the root $ \phi $. The \emph{multiplicity} of $ \phi $ is defined to be $ \dim \mathfrak{g}_\phi $, which is always finite.

The set $ \Delta $ of all roots of $ \mathfrak{g} $ is partitioned into the set $ \Delta_+ $ of \emph{positive roots}, consisting of the roots that are non-negative integral linear combinations of the simple roots, and the set $ \Delta_- $ of \emph{negative roots}, satisfying $ \Delta_- = - \Delta_+ $. The set $ \Pi $ of simple roots is partitioned into the set $ \Pi^\text{re} $ of \emph{real} roots $ \alpha_i $ such that $ a_{ii} > 0 $, and the set $ \Pi^\text{im} $ of \emph{imaginary} roots. 
There is a $ \mathbb{Z}^I $-gradation on $ \mathfrak{g}_0 $ and $ \mathfrak{g} $ given by 
\begin{equation*}
\deg e_i = (\delta_{ij})_{j \in I} , \quad \deg f_i = (-\delta_{ij})_{j \in I} , \quad \deg h_i = (0)_{j \in I} .
\end{equation*} 
The involution $ \omega $ on $ \mathfrak{g} $ given by
\begin{equation*}
\omega e_i = f_i, \ \omega f_i = e_i, \ \omega h_i = -h_i, \text{ for all } i \in I ,
\end{equation*}
negates the degree of each element and maps $ \mathfrak{g}_\phi $ onto $ \mathfrak{g}_{-\phi} $.

\section{Decomposition of Borcherds algebras with a zero block and all other entries negative}
\label{section:general decomposition}
In this section, we consider Borcherds algebras with no real roots and the property that any simple root that is orthogonal to another simple root is orthogonal to all such simple roots, including itself. In other words, we consider Borcherds algebras with Cartan matrices of the form 
\begin{equation}
\left(
\begin{blockarray}{c|c}
0 & <0 \\
\cline{1-2}
<0 & <0 \vspace{-1em}
\end{blockarray} 
\right).
\end{equation} 
We decompose these Borcherds algebras into free, Heisenberg and abelian subalgebras, which will provide an efficient calculation of twisted denominator formulas of the non-Fricke Monstrous Lie algebras in the final section. The results in this section are analogs of those in \cite{JurJPAA}. Though our result can be obtained as a corollary of Theorem 3.19 \cite{jurisich1996exposition}, proved by use of Lie algebra homology, we prove it here explicitly and directly. Here, the generators of the free and Heisenberg components will be explicitly seen. In the final section, the Monstrous Lie algebras associated to non-Fricke elements will be a special case of these results. 

We first state a result from \cite{bourbaki}. In what follows, $ L $ denotes the free Lie algebra functor on sets.
\begin{lemma}[The elimination theorem ] Let $ X $ be a set, $ S $ a subset of $ X $ and $ T $ the set of sequences $ (s_1, \dots, s_n, x) $ with $ n \geq 0 $, $ s_1, \dots, s_n $ in $ S $ and $ x \in X\setdif S $. Then the following statements are true.
\begin{enumerate}
\item[(a)] The Lie algebra $ L(X) $ is, as a vector space, the direct sum of the subalgebra $ L(S) $ of $ L(X) $ and the ideal $ \mathfrak{a} $ of $ L(X) $ generated by $ X \setdif S $.
\item[(b)] There exists a Lie algebra isomorphism 
\begin{equation*}
 \phi : L(T) \to \mathfrak{a} , \quad (s_1, \dots, s_n, x) \mapsto \ad(s_1) \cdots \ad(s_n) x .
\end{equation*}
\end{enumerate}
\end{lemma}

We also prove a lemma about how surjective set maps lift to Lie algebra maps of free Lie algebras.

\begin{lemma} \label{lemma:surj map}
Let $ f: X \to Y $ be a surjective set map and consider its image $ L(f) : L(X) \to L(Y) $ under the free Lie algebra functor. Then the kernel of $ L(f) $ is the ideal $ \mathfrak{i} $ in $ L(X) $ generated by $ \{ x_1 - x_2 : x_1, x_2 \in X \text{ with }f(x_1) = f(x_2) \} $. Furthermore, $ L(Y) \cong L(X)/ \mathfrak{i} $ as Lie algbebras.
\end{lemma}

\begin{proof}
By the universal property of free Lie algebras, we have 
\begin{equation}
\begin{tikzcd}
X \ar[d,"f",swap] \ar[r,"i_X"]  & L(X) \ar[d,dashed,"L(
f) = \widetilde{i_Y \circ f}"]\\ 
Y \ar[r,"i_Y",swap] & L(Y) .
\end{tikzcd}
\end{equation}
Since $ i_Y \circ f $ maps $ X $ onto the the set $ Y \subseteq L(Y) $ of generators of $ L(Y) $, we know that $ L(f) $ is surjective. Hence $ L(Y) \cong L(X)/\ker(L(f)) $, and we are left to study the kernel of $ L(f) $. 

Let $ \psi : Y \to L(X)/\mathfrak{i} $ be the set map that sends $ y \in Y $ to $ f^{-1}(y) + \mathfrak{i} \in L(X)/\mathfrak{i} $. Indeed, this is well-defined since $ f $ is surjective and for all $ x_1, x_2 \in X $ we have $ x_1 + \mathfrak{i} = x_2 + \mathfrak{i} $ when $ f(x_1) = f(x_2) $. Using the map $ \psi $ with the universal property of free Lie algebras
\begin{equation}
\begin{tikzcd}
Y \ar[r,"i_Y"] \ar[rd,"\psi",swap] & L(Y) \ar[d,dashed,"\widetilde{\psi}"]\\ 
& L(X)/\mathfrak{i} ,
\end{tikzcd}
\end{equation}
we obtain the Lie algebra map $ \widetilde{\psi} $. The ideal $ \mathfrak{i} $ is contained in the kernel of $ L(f) $, since the generators $ x_1 - x_2 \in \ker(L(f)) $ when $ f(x_1) = f(x_2) $. Hence, by the universal property of quotients 
\begin{equation}
\begin{tikzcd}
L(X) \ar[r,"\pi"] \ar[rd,"L(f)",swap] & L(X)/ \mathfrak{i} \ar[d,dashed,"\overline{L(f)}"]\\ 
& L(Y) ,
\end{tikzcd}
\end{equation}
we obtain the Lie algebra map $ \overline{L(f)} $. We then have the commutative diagram
\begin{equation}
\begin{tikzcd}
X \ar[r,"i_X"] \ar[dd,"f",swap] &L(X) \ar[d,"\pi"] \ar[dr,"L(f)",near end] \ar[drr,"\pi",bend left=15] \\
&L(X)/\mathfrak{i} \ar[r, "\overline{L(f)}",swap,near start] &L(Y) \ar[r,"\widetilde{\psi}",swap,near start] &L(X)/\mathfrak{i} .\\
Y \ar[ru,"\psi"] \ar[rru,"i_Y",bend right=15,shift left=1pt] \ar[rrru,"\psi",bend right=10, shift right=2pt,swap,near end]
\end{tikzcd}
\end{equation}
The top-right triangle commutes by the commutativity of the rest of the diagram together with the universal property of free Lie algebras. By the universal property of quotients, we have $ \widetilde{\psi} \circ \overline{L(f)} = \id_{L(X)/ \mathfrak{i}} $. So $ \overline{L(f)} $ is injective, hence $ \mathfrak{i} = \ker(L(f)) $.
\end{proof}

Using the previous lemmas, we prove the following theorem about Borcherds algebras defined from Cartan matrices with a diagonal block of zeros and all other entries being negative. This zero block will correspond to the sub-index set $ J $ in the theorem.

\begin{theorem} \label{theorem:decomposition}
Let $ A = (a_{ij})_{i,j \in I} $ be a symmetrized Cartan matrix with $ a_{ii} \leq 0 $ for all $ i \in I $. Let 
\begin{equation}
J = \{j \in I : a_{jj} = 0\}
\end{equation}
and assume that $ a_{ij} = 0 $ if and only if $ i,j \in J $. Then, given any ordering $ \leq $ on the indexing set $ J $, there are free subalgebras $ \mathfrak{u}^+ $ and $ \mathfrak{u}^- $ of $ \mathfrak{g} $ with generating sets
\begin{align}
U^+ &= \{ \ad(e_{j_1})^{\ell_1} \cdots  \ad(e_{j_n})^{\ell_n} e_i :  n \geq 0, \ j_1 < \dots < j_n \in J, \ \ell_k > 0, \ i \in I \setdif J \}, \\
U^- &= \{ \ad(f_{j_1})^{\ell_1} \cdots  \ad(f_{j_n})^{\ell_n} f_i :  n \geq 0, \ j_1 < \dots < j_n \in J, \ \ell_k > 0, \ i \in I \setdif J \}, 
\end{align}
respectively. Furthermore,
\begin{equation}
\mathfrak{g} = \mathfrak{g}(A) = \mathfrak{u}^- \oplus (\mathfrak{g}_J + \mathfrak{h}) \oplus \mathfrak{u}^+,
\end{equation}
where $ \mathfrak{g}_J $ is the Borcherds algebra associated to the matrix $ (a_{ij})_{i,j \in J} $.
\end{theorem}

\begin{proof}
We start with the triangular decompositions
\begin{equation*}
\mathfrak{g}_0 = \mathfrak{n}_0^- \oplus \mathfrak{h}_0 \oplus \mathfrak{n}_0^+ \qquad \text{and} \qquad \mathfrak{g} = \mathfrak{n}^- \oplus \mathfrak{h} \oplus \mathfrak{n}^+,
\end{equation*}
where $ \mathfrak{n}^{\pm} \cong \mathfrak{n}_0^{\pm}/ \mathfrak{k}_0^{\pm} $ and $ \mathfrak{h} \cong \mathfrak{h}_0 $. Using the elimination theorem on $ \mathfrak{n}_0^+ = L(X) $ with $ X = \{e_i: i \in I \} $ and $ S = \{e_j: j \in J \} $, we get
\begin{equation*}
\mathfrak{n}_0^+ = L(S) \oplus \mathfrak{a}
\end{equation*}
with $ \mathfrak{a} $ freely generated by
\begin{equation*}
\ad(e_{j_1}) \cdots \ad(e_{j_n}) e_i , \quad j_1, \dots, j_n \in J, \ i \in I \setdif J.
\end{equation*}

The subalgebra $ L(S) $, and  the ideals $ \mathfrak{a} $ and $ \mathfrak{k}_0^+ $ are graded with respect to the $ \mathbb{Z}^I $-gradation, since their generators are homogeneous. By the assumptions for the matrix $ A $, we know that $ \mathfrak{k}_0^+ $ is generated by $ \{[e_i,e_j] : i,j \in J \} $. Hence,
\begin{equation}
\mathfrak{k}_0^+ \cap L(S) =  \coprod_{\substack{ \mathbf{n} \in \mathbb{Z}^I \text{ s.t.} \\ n_i = 0 \text{ for all } i \in I \setdif J} } \mathfrak{k}_0^+(\mathbf{n}) 
\quad \text{and} \quad 
\mathfrak{k}_0^+ \cap \mathfrak{a} =  \coprod_{\substack{ \mathbf{n} \in \mathbb{Z}^I \text{ s.t.} \\ n_i \neq 0 \text{ for some } i \in I\setdif J} } \mathfrak{k}_0^+(\mathbf{n}),
\end{equation}
where $ \mathfrak{k}_0^+(\mathbf{n}) $ are the subspaces of $ \mathfrak{k}_0^+ $ of homogeneous elements of degree $ \mathbf{n} \in \mathbb{Z}^I $. So we have 
\begin{align}
\mathfrak{k}_0^+ = (\mathfrak{k}_0^+ \cap L(S)) \oplus (\mathfrak{k}_0^+ \cap \mathfrak{a}), 
\end{align}
and we can study the direct summands of
\begin{equation}
\mathfrak{n}^+ \cong L(S)/(\mathfrak{k}_0^+ \cap L(S)) \oplus \mathfrak{a}/ (\mathfrak{k}_0^+ \cap \mathfrak{a} ).
\end{equation}
Since $ \mathfrak{k}_0^+ \cap L(S) $ is the ideal in $ L(S) $ generated by $ \{ [e_i,e_j] : i,j \in J \} $, we have 
\begin{align}
 L(S)/(\mathfrak{k}_0^+ \cap L(S)) \cong \mathfrak{g}_J^+.
\end{align}

Choose and fix any ordering $ \leq $ on $ J $. Consider the set
\begin{align}
T' =  \{ (e_{j_1}, \dots, e_{j_n}, e_i) : n \geq 0, \ j_1 \leq \dots \leq j_n \in J, \ i \in I \setdif J \}
\end{align}
of non-decreasing sequences in $ \{e_j : j \in J \} $ ending with an element in $ \{e_i : i \in I \setdif J\} $. Consider the surjective set map $ \phi: T \to T' $ that sends sequences in $ T $ to their ordered version in $ T' $. The free Lie algebra functor $ L $ sends $ \phi $ to the surjective Lie algebra map $ L(\phi): L(T) \to L(T') $.
By Lemma \ref{lemma:surj map}, we know that $ \ker(L(\phi)) $ is the ideal $ \langle R \rangle $ of $ L(T) $ generated by the elements 
\begin{equation*}
R = \{t_1 - t_2 : t_1, t_2 \in T \text{ with } \phi(t_1) = \phi(t_2) \} .
\end{equation*}
After identifying $ L(T) $ with $ \mathfrak{a} $ via the elimination theorem, and repeatedly using
\begin{equation*}
 \ad(e_{j_1}) \ad(e_{j_2}) x = \ad(e_{j_2}) \ad(e_{j_1}) x + [[e_{j_1},e_{j_2}],x], \quad \text{for } {j_1},{j_2} \in J, \ x \in U^+, 
\end{equation*}
we can express the elements in $ R $ as elements in $ \mathfrak{k}_0^+ $. Hence, $ \langle R \rangle \subseteq \mathfrak{k}_0^+ \cap \mathfrak{a} $.

For the converse, we consider an element in $ b \in L(S) $ acting on an element in $ a \in \mathfrak{a} $ by the adjoint action. By repeated use of
\begin{equation*}
\ad(b)[x,y] = [\ad(b)x,y] + [x,\ad(b)y], \quad x,y \in \mathfrak{a} ,
\end{equation*}
we see that $ \ad(b)a $ is in the ideal of $ \mathfrak{a} $ generated by $ \ad(b)U^+ $. An element in $ \ad(\mathfrak{k}_0^+ \cap L(S))U^+ $ is a linear combination of elements of the form
\begin{equation*}
\ad(e_{j_1}) \cdots \ad([e_{j_k}, e_{j_{k+1}}]) \cdots \ad(e_{j_n}) x, \quad \text{for } j_1, \dots, j_n \in J, \ x \in U^+ .
\end{equation*}
These elements are seen to be in $ R $ after using 
\begin{equation*}
\ad([e_{j_k}, e_{j_{k+1}}]) = \ad(e_{j_k})\ad(e_{j_{k+1}}) - \ad(e_{j_{k+1}}) \ad(e_{j_k}).
\end{equation*}
Similar arguments show that the $ \ad(L(S))\langle R \rangle \subseteq \langle R \rangle $, hence $ \ad(\mathfrak{n}_0^+)\langle R \rangle \subseteq \langle R \rangle $.
Since $ \ad(\mathfrak{k}_0^+ \cap L(S))U^+ $ generates $ \mathfrak{k}_0^+ \cap \mathfrak{a} $, we have $ \mathfrak{k}_0^+ \cap \mathfrak{a}  \subseteq \langle R \rangle $. Hence, $ \ker(L(\phi)) =  \langle R \rangle 
= \mathfrak{k}_0^+ \cap \mathfrak{a} $, giving us
\begin{equation}\label{eq:iso}
\mathfrak{a}/(\mathfrak{k}_0^+ \cap \mathfrak{a}) \cong L(T)/\ker(L(\phi)) \cong L(T')
\end{equation}
as Lie algebras. We use $ T' \cong U^+ $ under $ (e_{j_1}, \dots, e_{j_n},e_i) \mapsto \ad(e_{j_1}) \cdots \ad(e_{j_n}) e_i $ to obtain $ L(T') \cong L(U^+) = \mathfrak{u}^+ $. Hence, so far, we have
\begin{align}
\mathfrak{n}^+ \cong \mathfrak{g}_J^+ \oplus \mathfrak{u}^+. 
\end{align} 
Using the involution $ \omega $, we similarly have
\begin{align}
\mathfrak{n}^- \cong \mathfrak{g}_J^- \oplus \mathfrak{u}^-.
\end{align}
Thus, we have the decomposition
\begin{equation}
\begin{aligned}
\mathfrak{g} &= \mathfrak{n}^- \oplus \mathfrak{h} \oplus \mathfrak{n}^+ = (\mathfrak{g}_J^- \oplus \mathfrak{u}^-) \oplus \mathfrak{h} \oplus (\mathfrak{g}_J^+ \oplus \mathfrak{u}^+) \\
&= \mathfrak{u}^- \oplus (\mathfrak{g}_J^- \oplus \mathfrak{g}_J^+ \oplus \mathfrak{h}) \oplus \mathfrak{u}^+ = \mathfrak{u}^- \oplus (\mathfrak{g}_J + \mathfrak{h}) \oplus \mathfrak{u}^+,
\end{aligned}
\end{equation}
as desired.
\end{proof}

We can determine the root system for $ \mathfrak{g} $ in the setting of the previous theorem. Since $ [e_i, e_j] = 0 $ when $ a_{ij} = 0 $, the elements in $ \mathfrak{g}_J^+ $ are linear combinations of $ \{e_j : j \in J\} $. Hence, the subalgebra $ \mathfrak{g}_J^+ $ contributes exactly the roots $ \{\alpha_j : j \in J \} $. The elements in 
\begin{align*}
U^+ &= \{ \ad(e_{j_1})^{\ell_1} \cdots  \ad(e_{j_n})^{\ell_n} e_i : n \geq 0, \ j_1 < \dots < j_n \in J, \ \ell_k > 0, \ i \in I \setdif J \}, 
\end{align*}
contribute exactly the roots
\begin{align*}
\{ \alpha_{j_1} + \dots + \alpha_{j_n} + \alpha_i :  n \geq 0, \ j_1 ,\dots , j_n \in J, \ i \in I \setdif J \} .
\end{align*}
Since $ U^+ $ freely generates $ \mathfrak{u}^+ $, we see that $ \mathfrak{u}^+ $ contributes exactly the roots
\begin{multline*}
\{ \alpha_{j_1} + \dots + \alpha_{j_n} + \alpha_{i_1} + \dots + \alpha_{i_p} :  n \geq 0, \ p > 0,\  j_k \in J,\ i_k \in I \setdif J \} \\  
\setdif \{ t \alpha_i : t \geq 2, \ i \in I \setdif J \}.
\end{multline*}
Applying the involution $ \omega $, we see that the subalgebras $ \mathfrak{g}_J^- $ and $ \mathfrak{u}^- $ contribute the negative of the roots above. 

Hence, we have the following result.
\begin{corollary}\label{theorem:positive roots}
In the setting of Theorem \ref{theorem:decomposition}, $ \mathfrak{g} $ has the set of positive roots
\begin{equation}
\begin{aligned}
\Delta_+ &= \{\alpha_j : j \in J \} \,  \sqcup  \\
& \ \{ \alpha_{j_1} + \dots + \alpha_{j_n} + \alpha_{i_1} + \dots + \alpha_{i_p} :  n \geq 0, \ p > 0,\  j_k \in J,\ i_k \in I \setdif J \} \\   
&\  \setdif \{ t \alpha_i : t \geq 2, \ i \in I \setdif J \}.
\end{aligned}
\end{equation}
\end{corollary}

In the final section, we will be interested in quotienting $ \mathfrak{g} $ by its center $ \mathfrak{z} $. In certain cases, as in our main application, $ \mathfrak{g}/\mathfrak{z} $ will have a finite-dimensional Cartan subalgebra. 

Recall from Section \ref{section:GKMs} that we defined roots of $ \mathfrak{g} $ to be elements in $ (\mathfrak{h}^e)^* $ so that the simple roots are linearly independent. We can define the \emph{simple roots} $ \projh{\alpha}_i $, $ i \in I $, as elements of $ \mathfrak{h}^* $ satisfying
\begin{equation}
[h, e_i] = \projh{\alpha}_i(h) e_i, \text{ for all } h \in \mathfrak{h} .
\end{equation}
Analogously, we can define a \emph{root} to be $ \phi \neq 0 \in \mathfrak{h}^* $ with a non-zero \emph{root space}
\begin{equation}
\mathfrak{g}_\phi = \{ x \in \mathfrak{g} : [h, x] = \phi(h) x, \text{ for all } h \in \mathfrak{h} \}.
\end{equation} 
It is possible that the simple roots are linearly dependent as elements of $ \mathfrak{h}^* $. Moreover, it is possible that $ \projh{\alpha}_i  = \projh{\alpha}_j $ for some $ i \neq j $. In this case, the dimension of a root space for this simple root is greater than one. 

Since $ \mathfrak{h}^e = \mathfrak{h} \oplus \mathfrak{d} $, there is a linear map 
\begin{equation}
\projh{\,\cdot\,} : (\mathfrak{h}^e)^* \to \mathfrak{h}^*, \quad \projh{\phi}: h_i \mapsto \phi(h_i),
\end{equation}
which restricts the domain of the functionals. This gives a precise way to interpret the roots of $ \mathfrak{g} $ as elements in $ \mathfrak{h}^* $ instead of elements in $(\mathfrak{h}^e)^* $.

Let $ B \subseteq I $ be such that $ \projh{\alpha}_i \neq \projh{\alpha}_j $, for all $ i \neq j $ in $ B $, and such that $ \{ \projh{\alpha}_i\}_{i \in B} $ is a basis for the subspace $ \Span \{ \projh{\alpha}_i \}_{i \in I} $ of $ \mathfrak{h}^* $. 
\begin{lemma} \label{lemma:center}
We have the decomposition 
\begin{equation}
\mathfrak{h} = \Span \{ h_i \}_{i \in B} \oplus \mathfrak{z},
\end{equation}
where $ \mathfrak{z} $ is the center of $ \mathfrak{g} $.
\end{lemma}

\begin{proof} 
It is known that $ \mathfrak{z} \subseteq \mathfrak{h} $, so we show that $ \mathfrak{h} = \Span \{ h_i \}_{i \in B} + \mathfrak{z} $ and $  \Span \{ h_i \}_{i \in B} \cap \mathfrak{z} = 0 $. Since $ \{\projh{\alpha}_i\}_{i \in B} $ is a basis for $ \Span \{ \projh{\alpha}_i \}_{i \in I} $, we have $ (c_i^{(j)})_{i \in B} $ such that $ \projh{\alpha}_j = \sum_{i \in B} c_i^{(j)} \projh{\alpha}_i $, for each $ j \in I \setdif B $. Define
\begin{align*}
\widetilde{h}_j = h_j - \sum_{i \in B} c_i^{(j)} h_i, \text{ for all } j \in I \setdif B.
\end{align*}
Since 
\begin{align*}
[\widetilde{h}_j, e_k] &= \left( \projh{\alpha}_k(h_j) - \sum_{i \in B} c_i^{(j)} \projh{\alpha}_k(h_i)\right) e_k = \left( \projh{\alpha}_j(h_k) - \sum_{i \in B} c_i^{(j)} \projh{\alpha}_i(h_k)\right) e_k = 0,
\end{align*} 
and similarly for $ f_k $, we see that $ \{ \widetilde{h}_j \}_{j\in I \setdif B} $ is contained in $ \mathfrak{z} $.
 Hence,
\begin{equation*}
\mathfrak{h} \subseteq \Span\{h_i\}_{i \in B} +\Span\{ \widetilde{h}_j \}_{j\in I \setdif B} \subseteq \Span\{h_i\}_{i \in B} + \mathfrak{z} \subseteq \mathfrak{h}.
\end{equation*}  
Furthermore, 
if $ h = \sum_{i \in B} c_i h_i $ is an arbitrary element in $  \Span\{h_i\}_{i \in B}  \cap \mathfrak{z} $, then 
\begin{align*}
0 = [ h, e_j] = \sum_{i \in B} c_i  \projh{\alpha}_j( h_i) e_j
\end{align*}
implies 
\begin{align*}
\sum_{i \in B} c_i \projh{\alpha}_i(h_j) = \sum_{i \in B} c_i  \projh{\alpha}_j( h_i)  = 0, \text{ for all } j \in I.
\end{align*}
Hence, $ \sum_{i \in B} c_i \projh{\alpha}_i = 0 $ implying $ c_i = 0 $ for all $ i \in B $, by the linear independence of $ \{\projh{\alpha}_i\}_{i \in B} $. So $ h = 0 $.
\end{proof}
Using Theorem \ref{theorem:decomposition} and Lemma \ref{lemma:center}, we have the following proposition.

\begin{proposition} \label{theorem:quotient structure}
Let $ A $ be as in the setting of Theorem \ref{theorem:decomposition}. Then the quotient of $ \mathfrak{g}(A) $ by its center has the decomposition
\begin{equation}
\mathfrak{g}(A)/\mathfrak{z} = (\mathfrak{u}^- \oplus \mathfrak{g}_J^- ) \oplus \Span \{h_i\}_{i \in B} \oplus (\mathfrak{u}^+ \oplus \mathfrak{g}_J^+ ),
\end{equation}
where $ B \subseteq I $ is such that $ \projh{\alpha}_i \neq \projh{\alpha}_j $, for all $ i \neq j $ in $ B $, and $ \{ \projh{\alpha}_i \}_{i \in B} $ is a basis for $ \Span\{ \projh{\alpha}_i \}_{i \in I} \subseteq \mathfrak{h}^* $. \qed
\end{proposition}

Since $ \mathfrak{z} $ is a graded ideal with respect to the $ \mathbb{Z}^I $-gradation, the quotient $ \mathfrak{g}(A)/\mathfrak{z} $ inherits this grading. The degree zero component $ \overline{\mathfrak{h}} = \Span \{h_i\}_{i \in B} $ is the Cartan subalgebra of $ \mathfrak{g}/\mathfrak{z} $, which is finite-dimensional if and only if $ B $ is finite. The roots and simple roots of $ \mathfrak{g}(A)/\mathfrak{z} $ can now be defined as elements of $ \overline{\mathfrak{h}}^* $. It is possible that the simple roots as elements of $ \overline{\mathfrak{h}}^* $ are linearly dependent and that the dimension of a root space for a simple root is greater than one. Since roots in $ (\mathfrak{h}^e)^* $ vanish on $ \mathfrak{z} $, there is a well-defined linear map 
\begin{equation}\label{eq:hat proj}
\projq{\,\cdot\,} : (\mathfrak{h}^e)^* \to \overline{\mathfrak{h}}^*, \quad \projq{\phi}: h_i + \mathfrak{z} \mapsto \phi(h_i),
\end{equation}  
which interprets the roots of $ \mathfrak{g} $ as elements of $ \overline{\mathfrak{h}}^* $. 
Since $ \{h_i\}_{i\in B} \cup \{ \widetilde{h}_j \}_{j\in I \setdif B} $ generates $ \mathfrak{h} $, a linear combination $ \sum_{i \in B} c_i \projq{\alpha}_i $ that is zero on $ \Span\{h_i\}_{i\in B} $ implies that $ \sum_{i \in B} c_i \projh{\alpha}_i $ is zero on all of $ \mathfrak{h} $. So, the linear independence of $ \{\projh{\alpha}_i\}_{i\in B} $ gives the linear independence of $ \{ \projq{\alpha}_i \}_{i \in B} $, which is hence basis for $ \Span\{\projq{\alpha}_i \}_{i \in I} $.

\section{Application to non-Fricke Monstrous Lie algebras} \label{section:monstrous Lie algebras}

The Monstrous Lie algebras $ \mathfrak{m}_g $, as defined by Carnahan in \cite{CarThesis}, \cite{CarDuke} and \cite{CarIV}, are certain Lie algebras associated to each element $ g $ in the Monster group $ \M $. They can be defined by an explicit Cartan matrix $ A_g $ as $ \mathfrak{m}_g = \mathfrak{g}(A_g)/\mathfrak{z} $, where $ \mathfrak{z} $ is the center of $ \mathfrak{g}(A_g) $ as outlined in \cite{CarDuke} (with notation $ W_g $). In \cite{CarIV}, Caranahan showed that this Lie algebra is isomorphic to a Lie algebra obtained from an abelian intertwining algebra of $ g^i $-twisted $ V^\natural $-modules, $ i = 0, \dots, o(g)-1 $, by use of a certain quantization functor. The matrix $ A_g $ satisfies the assumptions in Theorem \ref{theorem:decomposition} when $ g $ is of non-Fricke type. This is our main application of the results from the previous section.

To each element $ g \in \M $, we have the \emph{McKay--Thompson} series
\begin{align*}
T_g(\tau) = \sum_{n = 0}^\infty \Tr\left( g|_{V^\natural_n} \right) q^{n-1}, \quad q = e^{2\pi i \tau} .
\end{align*} 
An element $ g \in \M $ is called \emph{Fricke} if $ T_g(\tau) $ is invariant under the level $ N $ Fricke involution $ \tau \mapsto -1/N\tau $ for some $ N \in \mathbb{Z}_{>0} $. Otherwise, $ g $ is called \emph{non-Fricke}. When $ g $ is Fricke, $ \mathfrak{m}_g $ has structure similar to the Monster Lie algebra $ \mathfrak{m} = \mathfrak{m}_1 $ as discussed in \cite{CarFricke}.  We will present the Cartan matrix for the case that $ g $ is non-Fricke. 

We define the sequence $ \{c(m,0)\}_{m=1}^\infty $ of non-negative integers by 
\begin{equation}\label{eq:c(m,0) defn}
T_g(\tau) = q^{-1} \prod_{m>0} (1- q^m)^{c(m,0)} .
\end{equation}
This sequence can be shown to exist and computed explicitly for each $ g $ using a known $ \eta $-product expression
\begin{equation}
 T_g(\tau) = \eta(a_1 \tau)^{b_1} \cdots \eta(a_s\tau)^{b_s}, 
\end{equation}
where $ \eta(\tau) = q^{1/24} \prod_{n=1}^\infty (1 - q^n) $ is the Dedekind $ \eta $-function and $ a_1, b_1, \dots, a_s, b_s $ are certain integers. These integers can be found in Tables 1 and 2 of \cite{TuiteMon}. By expanding this product, we can see the coefficient of $ q^m $ is the sum of $ -c(m,0) $ and binomial coefficients of $ c(n,0) $ for $ n < m $. Hence, such a sequence of non-negative integers $ \{c(m,0)\}_{m=1}^\infty $ can be inductively shown to be unique.

Using the functional equation $ \eta(-1/\tau) = \sqrt{-i\tau} \eta(\tau) $ on the $ \eta $-product, we find the generating function
\begin{equation}\label{eq:c(1,n/N) defn}
T_g(-1/\tau) = \sum_{n>0} c(1, n/N) q^{n/N} = (a_1^{b_1} \cdots a_s^{b_s})^{-1/2} q^{1/N} \prod_{m > 0} (1 - q^{m/N})^{-c(m,0)} 
\end{equation}
for the sequence $ \{c(1,n/N) \}_{n=1}^\infty $ of non-negative integers, for some $ N \in \mathbb{Z}_{>0} $. Furthermore, $ c(1, 1/N) = (a_1^{b_1} \cdots a_s^{b_s})^{-1/2} $ is non-zero.  Since the  coefficient of $ q^{1/N} $ is non-zero and $ T_g(-1/\tau) $ converges in some open region, the value $ N $ and the sequence $ \{c(1,n/N) \}_{n=1}^\infty $ of non-negative integers is uniquely defined by $ T_g(-1/\tau) = \sum_{n>0} c(1, n/N) q^{n/N} $ when $ c(1, 1/N) $ is assumed to be non-zero. The computations and proof of existence of these sequences can be found in \cite{CCDHTT}. 

We fix a non-Fricke element $ g \in \M $, and consider the sequences  $ \{c(m,0)\}_{m=1}^\infty $  and $ \{c(1,n/N) \}_{n=1}^\infty $ that correspond to it. We define the index sets
\begin{align*}
I_- &= \{ (-m,j) \in \mathbb{Z} \times \mathbb{Z} : m > 0 , \ 1 \leq j \leq c(m,0) \},  \\ 
I_+ &= \{ (n,k) \in \mathbb{Z} \times \mathbb{Z} : n > 0 , \ 1 \leq k \leq c(1,n/N) \}, 
\end{align*}
and $ I = I_- \cup I_+ $. The ordering for this index set is $ (m,j) \leq (n,k) $ if and only if $ m \leq n $ and $ j \leq k $. It is possible that $ c(m,0) = 0 $ or $ c(1,n/N) = 0 $ for some $ m,n > 0 $. In these cases, $ (-m,j) $ and $ (n,k) $ do not appear in the index set $ I $. It is always the case that $ c(m,0) \neq 0 $ for some $ m > 0 $ and $ c(1,1/N) \neq 0 $.

For each non-Fricke element $ g \in \M $, we define $ A_g $ to be the $ I \times I $ matrix with 
\begin{align*}
a_{(-m,j),(-n,k)} &= 0, && \text{when } (-m,j),(-n,k) \in I_-, \\
a_{(-m,j),(n,k)} &= -mn , && \text{when } (-m,j) \in I_- \text{ and } (n,k) \in I_+, \\
a_{(m,j),(-n,k)} &= -mn , && \text{when } (m,j) \in I_+ \text{ and } (-n,k) \in I_-, \\ 
a_{(m,j),(n,k)} &= -(m+n) , && \text{when } (m,j), (n,k) \in I_+ . 
\end{align*}
In other words, $ A_g $ is the matrix\\
\resizebox{1\textwidth}{!} 
{
$
\begin{blockarray}{(ccc|ccc)}
\ddots  &\vdots &\vdots
&\vdots&\vdots&\iddots\\
\cdots &0_{c(3,0)\times c(2,0)} &0_{c(3,0)\times c(1,0)}  
&-3_{c(3,0)\times c(1,1/N)}&-6_{c(3,0)\times c(1,2/N)} &\cdots\\
\cdots &0_{c(2,0)\times c(2,0)} &0_{c(2,0)\times c(1,0)}
&-2_{c(2,0)\times c(1,1/N)}&-4_{c(2,0)\times c(1,2/N)}&\cdots\\
\cdots  & 0_{c(1,0)\times c(2,0)} &   0_{c(1,0)\times c(1,0)}  
&-1_{c(1,0)\times c(1,1/N)}&-2_{c(1,0)\times c(1,2/N)} &\cdots\\
\cline{1-6}
\cdots  &-2_{c(1,1/N)\times c(2,0)} &  -1_{c(1,1/N)\times c(1,0)}
&-2_{c(1,1/N)\times c(1,1/N)}& -3_{c(1,1/N)\times c(1,2/N)} &\cdots\\
\cdots  &-4_{c(1,2/N)\times c(2,0)}&-2_{c(1,2/N)\times c(1,0)} 
&-3_{c(1,2/N)\times c(1,1/N)}&-4_{c(1,2/N)\times c(1,2/N)} &\cdots\\
\cdots &-6_{c(1,3/N)\times c(2,0)} &-3_{c(1,3/N)\times c(1,0)}  
&-4_{c(1,3/N)\times c(1,1/N)}&-5_{c(1,3/N)\times c(1,2/N)} &\cdots\\
\iddots &\vdots &\vdots  
&\vdots&\vdots &\ddots
\end{blockarray}\ ,$
} where $ r_{m \times n} $ denotes the $ m \times n $ matrix with all entries taking the values $ r $. The indexing for $ I_+ $ increases in the down and right directions from the center. Hence, $ I_- $ decreases in the up and left directions from the center. The matrix $ A_g $ is indeed a symmetrized Cartan matrix, since it is symmetric with non-negative entries.

We define the \emph{monstrous Lie algebra} $ \mathfrak{m}_g $ associated to $ g $ to be the Lie algebra 
\begin{equation*}
\mathfrak{m}_g = \mathfrak{g}(A_g) / \mathfrak{z},
\end{equation*}
where $ \mathfrak{z} $ is the center of $ \mathfrak{g}(A_g) $.

We see that $ A_g $ satisfies the assumptions in Theorem \ref{theorem:decomposition} with $ J = I_- $ and $ I \setdif J = I_+ $. Hence, we can directly apply this theorem to uncover the structure of $ \mathfrak{g}(A_g) $.

\begin{proposition} \label{theorem:g(A_g) structure}
There are free subalgebras $ \mathfrak{u}^+ $ and $ \mathfrak{u}^- $ of $ \mathfrak{g}(A_g) $ with generating sets
\begin{align}
U^+ &= \{ \ad(e_{j_1})^{\ell_1} \cdots  \ad(e_{j_n})^{\ell_n} e_i :  n \geq 0, \ j_1 < \dots < j_n \in I_-, \ \ell_k > 0, \ i \in I_+ \}, \\
U^- &= \{ \ad(f_{j_1})^{\ell_1} \cdots  \ad(f_{j_n})^{\ell_n} f_i :  n \geq 0, \ j_1 < \dots < j_n \in I_-, \ \ell_k > 0, \ i \in I_+ \}, 
\end{align}
respectively. Furthermore,
\begin{equation}
\mathfrak{g}(A_g) = \mathfrak{u}^- \oplus (\mathfrak{g}_{I_-} + \mathfrak{h}) \oplus \mathfrak{u}^+,
\end{equation}
where $ \mathfrak{g}_{I_-} $ is the Borcherds algebra associated to the matrix $ (a_{ij})_{i,j \in I_-} $. \qed
\end{proposition}

We now proceed to find a basis of $ \Span\{\projh{\alpha}_i\}_{i\in I} \subseteq \mathfrak{h}^* $. Recall that the simple roots $ \{\projh{\alpha}_i\}_{i\in I} $ are possibly linearly dependent since they are elements in $ \mathfrak{h}^* $, and not the linearly independent roots $ \{\alpha_i\}_{i\in I} $ in $ (\mathfrak{h}^e)^* $.

\begin{lemma} \label{lemma:basis}
For some $ \m > 0 $, $ \{\projh{\alpha}_{(-\m,1)}, \projh{\alpha}_{(1,1)} \} $ is a basis for $ \Span\{\projh{\alpha}_i\}_{i\in I} \subseteq \mathfrak{h}^* $.
\end{lemma}

\begin{proof} Consider the case when $ c(1,0) \neq 0 $. To make the computations more concrete, we identify the rows of $ A_g $ with $ \Span\{\projh{\alpha}_i\}_{i\in I} $ by corresponding the $ i^{\text{th}}$ row with $ \projh{\alpha}_i $.

Firstly, $ \{ r_{(-1,1)}, r_{(1,1)} \} $ is linearly independent because $ r_{(-1,1)} $ has zero entries where $ r_{(1,1)} $ does not. Next we show that $ \{ r_{(-1,1)}, r_{(1,1)} \} $ spans the row space by showing that every row $ r_i $, $ i \in I $, is a linear combination of these two rows. \\
Consider the case when $ i = (-m,k) \in I_- $. Let $ j = (n,q) \in I $. Then
\begin{align*}
(r_i)_j = a_{i,j} = \begin{cases}
0, & j \in I_-, \\
-mn, &j \in I_+ 
\end{cases} 
= \begin{cases}
0, & j \in I_-, \\
m(-1 \cdot n), &j \in I_+ 
\end{cases}
= m a_{(-1,1),j}.
\end{align*}
So $ r_i = m r_{(-1,1)} $. Consider the case when $ i = (m,k) \in I_+ $. Let $ j = (n,q) \in I $. Then
\begin{align*}
(r_i)_j &= a_{i,j} = \begin{cases}
mn, & j \in I_-, \\
-(m+n), &j \in I_+ 
\end{cases} 
= \begin{cases}
m(1\cdot n)+0, & j \in I_-, \\
(1-m)(-1 \cdot n)-m(1+n) , &j \in I_+ 
\end{cases} \\
&=  (1-m)a_{(-1,1),j} + m a_{(1,1),j}
\end{align*}
So  $ r_i =  (1-m) r_{(-1,1)} + m r_{(1,1)} $. Hence, every row of $ A $ is a linear combination of the rows $ r_{(-1,1)} $ and $ r_{(1,1)} $.

In the case that $ c(1,0) = 0 $, we pick $ \m $ to be the minimal positive integer such that $ c(\m,0) \neq 0 $. Since the $ (-\m,1)^\text{th} $ row is a multiple of the (non-existent) $ (-1,1)^\text{th} $ row, the argument proceeds similarly as above.
\end{proof}

In the notation of Proposition \ref{theorem:quotient structure}, this lemma gives $ B = \{ (-\m,1), (1,1) \} $ where $ \m > 0$ is such that $ c(\m,0) \neq 0 $. We will assume that $ \m $ is minimal by convention. 

\begin{theorem}\label{theorem:m_g structure}
The non-Fricke Monstrous Lie algebras can be expressed as
\begin{equation}
\mathfrak{m}_g = \mathfrak{u}^- \oplus \heis \oplus \mathbb{C} h_{(1,1)}  \oplus \mathfrak{u}^+,
\end{equation}
where $ \mathfrak{u}^\pm $ are as in Proposition \ref{theorem:g(A_g) structure}, $ \heis $ is the infinite-dimensional Heisenberg Lie algebra spanned by $ \{e_i, f_i\}_{i \in I_-} $ and its central element $ h_{(-\ell,1)} $. 
\end{theorem}

\begin{proof}
By using Proposition \ref{theorem:quotient structure}, we obtain 
\begin{align*}
\mathfrak{g}/\mathfrak{z} &= (\mathfrak{u}^- \oplus \mathfrak{g}_{I_-}^- ) \oplus ( \mathbb{C} h_{(-\m,1)} \oplus \mathbb{C} h_{(1,1)} ) \oplus ( \mathfrak{g}_{I_-}^+ \oplus \mathfrak{u}^+ )\\
&= \mathfrak{u}^- \oplus (\heis \oplus \mathbb{C} h_{(1,1)} ) \oplus \mathfrak{u}^+,
\end{align*}
where $ \heis =  \mathfrak{g}_{I_-}^- \oplus \mathbb{C} h_{(-\m,1)} \oplus \mathfrak{g}_{I_-}^+ $. 
For each $ j = (-m,k) \in {I_-} $, in $ \mathfrak{g} $ we have
\begin{align*}
[e_j, f_j] = h_j = mh_{(-\m,1)} + h,
\end{align*}
where $ h \in \mathfrak{z} $. Furthermore, $ [h_{(-\m,1)},e_j] = [h_{(-\m,1)},f_j] = 0 $, making the derived subalgebra of $ \heis $ equal to its one-dimensional center. Hence, $ \heis $ is indeed a Heisenberg Lie algebra. 
\end{proof}

The infinite-dimensional Heisenberg subalgebra $ \heis $ of $ \mathfrak{m}_g $, when $ g $ is non-Fricke, comes from the infinite zero block in the Cartan matrix. This is analogous to the subalgebra $ \mathfrak{sl}_2 $ of $ \mathfrak{m}_g $ when $ g $ is Fricke, which comes from a single $ 2 $ in the Cartan matrix \cite{JurJPAA}, \cite{CarFricke}.

Consider the $ \heis $-modules
\begin{align}\label{eq:V+-}
\V^+ = \coprod_{i \in I_+} \mathcal{U}(\heis) \otimes_{\mathfrak{g}_{I_-}^- \oplus \mathbb{C}h_{-\m,1}} \mathbb{C} e_i \quad \text{and} \quad 
\V^- = \coprod_{i \in I_+} \mathcal{U}(\heis) \otimes_{\mathbb{C}h_{-\m,1} \oplus \mathfrak{g}_{I_-}^+} \mathbb{C} f_i,
\end{align}
where $\mathfrak{g}_{I_-}^- \oplus \mathbb{C}h_{-\m,1} $ acts on $ e_i $ via the adjoint action in $ \mathfrak{m}_g $, and $ \mathbb{C}h_{-\m,1} \oplus \mathfrak{g}_{I_-}^+ $ acts on $ f_i $ via the adjoint action in $ \mathfrak{m}_g $. By the PBW theorem, $ U^\pm $ is identifiable with a basis for $ \V^\pm $. Since the action of $ \heis $ on $ \V^\pm $ agrees with the adjoint action, we can identify $ \V^\pm $ with $ \Span U^\pm $. Hence, we have the immediate corollary. 

\begin{corollary}
The subalgebras $ \mathfrak{u}^\pm $ of $ \mathfrak{m}_g $ are freely generated by the vector spaces $ \V^\pm $.
\end{corollary}

Now consider the triangular decomposition
\begin{align*}
\mathfrak{m}_g = \left(\mathfrak{u}^- \oplus \mathfrak{g}_J^- \right) 
\oplus
\left( \mathbb{C} h_{(-\m,1)} \oplus \mathbb{C} h_{(1,1)} \right) 
\oplus 
\left( \mathfrak{u}^+\oplus ( \mathfrak{g}_J)^+ \right) = \mathfrak{g}^- \oplus \overline{\mathfrak{h}} \oplus \mathfrak{g}^+
\end{align*}
with a two-dimensional Cartan subalgebra $ \overline{\mathfrak{h}} = \mathbb{C} h_{(-\m,1)} \oplus \mathbb{C} h_{(1,1)} $. Recall that we use $ \overline{\mathfrak{h}} $ to denote the Cartan subalgebra of $\mathfrak{g} / \mathfrak{z} $ and $ \projq{\alpha}_i $ to denote the simple roots as elements in $ \overline{\mathfrak{h}}^* $.  By the general results in  Section \ref{section:general decomposition}, we know the two simple roots $ \projq{\alpha}_{(-\m,1)} $ and $ \projq{\alpha}_{(1,1)} $ form a basis for $ \overline{\mathfrak{h}}^* $, and the proof of Lemma \ref{lemma:basis} gives the remaining simple roots 
\begin{align*}
\projq{\alpha}_{(-m,j)} &= \frac{m}{\m} \projq{\alpha}_{(-\m,1)}, &\text{for all } (-m,j) \in I_- ,\\
\projq{\alpha}_{(n,k)} &= \frac{(1-n)}{\m}\projq{\alpha}_{(-\m,1)} + n \projq{\alpha}_{(1,1)}, &\text{for all } (n,k) \in I_+.
\end{align*}
Since $ \projq{\alpha}_{(-m,j)} $, for $ j = 1, \dots, c(m,0) $, are all the same root and no other roots are equal to it, we know that the multiplicity of $ \projq{\alpha}_{(-m,j)} $ is $ c(m,0) $. Similarly, the multiplicity of $ \projq{\alpha}_{(n,k)} $ is $ c(1,n/N) $. 

Applying Corollary \ref{theorem:positive roots}, we know that the positive roots are of the form
\begin{align*}
\projq{\alpha} &=  \frac{(m - n)}{\m} \projq{\alpha}_{(-\m,1)} +  n \projq{\alpha}_{(1,1)} ,   \ \text{for some } m \in \mathbb{Z}_{>0}, \ n \in \mathbb{Z}_{\geq 0} .
\end{align*}
Note that we have applied the projection \eqref{eq:hat proj}, so it is possible that $ t \projq{\alpha}_{(n,k)} $, for some $ t \geq 2 $, is a root if it is obtained by a linear combination with a root $ \projq{\alpha}_{(-m,j)} $ or if the multiplicity of $ \projq{\alpha}_{(n,k)} $ is greater than one. We expand the notation of $ c(m,0) $ and $ c(1, n/N) $ and denote the multiplicity of the root $ \frac{(m  - n)}{\m} \projq{\alpha}_{(-\m,1)} +  n \projq{\alpha}_{(1,1)} $ by $ c(m,n/N) $. Furthermore, we define $ c(m,n/N) = 0 $ if the functional $ \frac{(m  - n)}{\m} \projq{\alpha}_{(-\m,1)} +  n \projq{\alpha}_{(1,1)} $ is not a root.  All root multiplicities are finite since the multiplicities of the simple roots are finite.

In the following section, we will see in equation \eqref{eq:denomIDnon-fricke} that the $ c(m,n/N) $ appear in a generating function relation containing both the $ c(m,0) $ and $ c(1, n/N) $.

By plotting the root $ \frac{(m  - n)}{\m} \projq{\alpha}_{(-\m,1)} +  n \projq{\alpha}_{(1,1)} $ at $ (m,n/N) $ on the lattice $ \mathbb{Z} \times \frac{1}{N} \mathbb{Z} $, we can see how the two-dimensional Cartan subalgebra grades $ \mathfrak{m}_g $. Note that the lattice obtained when grading by weight is naturally $ \mathbb{Z} \times \mathbb{Z} $. However, we can simply scale one of the components to match the power of $ q $ in $ T_g(-1/\tau) $ as we will see in the next section. We can add $ \mathbb{C} h_{(1,1)} $ to $ \heis $, so that $ \heis \oplus \mathbb{C} h_{(1,1)} $ contains the Cartan subalgebra, and view $ \V^\pm $ as $ \heis \oplus \mathbb{C} h_{(1,1)} $-modules via the adjoint action. With this grading, we can see how $ \V^\pm $ are in fact direct sums of lowest/highest-weight $ \heis $-modules. In the root systems, the lowest-weight $ \heis $-modules are contained in the strings moving to the right from the roots $ \projq{\alpha}_{(n,k)} $. The subalgebra $ \mathfrak{u}^+ $ is then freely generated from direct sum of these lowest-weight $ \heis $-modules. On the horizontal axis lies only the positive part $ \mathfrak{g}_{I_-}^+ $ of $ \mathcal{H} $.

As an example, we look at the root systems of $ \mathfrak{m}_g $ when $ g $ is in the conjugacy class $ 2B $ or $ 4D $ (see Figures \ref{fig:2B} and \ref{fig:4D}). Since the negative roots are the negation of the positive roots, we only look at the positive half. The $ \eta $-product expressions $ T_{2B}(\tau) = \eta(\tau)^{24} \eta(2\tau)^{-24} $ and $ T_{4D} (\tau) = \eta(2 \tau)^{12} \eta(4 \tau)^{-12} $  can be used to compute the simple root multiplicities 
\begin{gather*}
T_{2B}(\tau) = q^{-1} \prod_{\text{odd } m > 0}(1-q^m)^{24} \quad \text{and}\quad T_{4D}(\tau) = q^{-1} \prod_{\substack{ m > 0 \\ m = 2 \text{ mod }4}}(1-q^m)^{12} , \\
T_{2B}(-1/\tau) = 2^{12}q^{1/2} \prod_{\text{odd } m > 0}(1-q^{m/2})^{-24}, \\
T_{4D}(-1/\tau) = 2^{6}q^{1/8} \prod_{\substack{ m > 0 \\ m = 2 \text{ mod }4}}(1-q^{m/8})^{-12},
\end{gather*}
as outlined in the beginning of this section. 
\begin{figure}[t]
\centering
{\footnotesize
\begin{tikzpicture}
	\begin{pgfonlayer}{nodelayer}
		\node [style=none] (0) at (0, 4.75) {};
		\node [style=none] (1) at (0, -0.5) {};
		\node [style=none] (2) at (7.75, 0) {};
		\node [style=none] (3) at (-0.5, 0) {};
		\node [style=root1] (4) at (1, 0) {};
		\node [style=root1] (6) at (7, 0) {};
		\node [style=root2] (12) at (1, 1) {};
		\node [style=root2] (13) at (1, 2) {};
		\node [style=root2] (14) at (1, 3) {};
		\node [style=root2] (15) at (1, 4) {};
		\node [style=root3] (21) at (3, 1) {};
		\node [style=root3] (22) at (5, 1) {};
		\node [style=root3] (24) at (3, 2) {};
		\node [style=root3] (25) at (5, 2) {};
		\node [style=root3] (27) at (3, 3) {};
		\node [style=root3] (28) at (5, 3) {};
		\node [style=root3] (30) at (3, 4) {};
		\node [style=root3] (31) at (5, 4) {};
		\node [style=none] (52) at (-0.5, 1) {$1/2$};
		\node [style=none] (53) at (-0.5, 2) {$2/2$};
		\node [style=none] (54) at (-0.5, 3) {$3/2$};
		\node [style=none] (55) at (-0.5, 4) {$4/2$};
		\node [style=none] (60) at (1, -0.5) {$1$};
		\node [style=none] (61) at (2, -0.5) {$2$};
		\node [style=none] (62) at (3, -0.5) {$3$};
		\node [style=none] (63) at (4, -0.5) {$4$};
		\node [style=none] (68) at (1, 0.5) {$\projq{\alpha}_{(-1,j)}$};
		\node [style=none] (70) at (7, 0.5) {$\projq{\alpha}_{(-7,j)}$};
		\node [style=none] (72) at (1, 1.5) {$\projq{\alpha}_{(1,k)}$};
		\node [style=none] (73) at (1, 2.5) {$\projq{\alpha}_{(2,k)}$};
		\node [style=none] (74) at (1, 3.5) {$\projq{\alpha}_{(3,k)}$};
		\node [style=none] (75) at (1, 4.5) {$\projq{\alpha}_{(4,k)}$};
		\node [style=none] (76) at (5, -0.5) {$5$};
		\node [style=none] (77) at (6, -0.5) {$6$};
		\node [style=none] (78) at (7, -0.5) {$7$};
		\node [style=root3] (80) at (7, 1) {};
		\node [style=root3] (81) at (7, 2) {};
		\node [style=root3] (82) at (7, 3) {};
		\node [style=root3] (83) at (7, 4) {};
		\node [style=root1] (84) at (3, 0) {};
		\node [style=root1] (85) at (5, 0) {};
		\node [style=none] (86) at (3, 0.5) {$\projq{\alpha}_{(-3,j)}$};
		\node [style=none] (87) at (5, 0.5) {$\projq{\alpha}_{(-5,j)}$};
		\node [style=root3] (88) at (2, 1) {};
		\node [style=root3] (89) at (4, 1) {};
		\node [style=root3] (90) at (6, 1) {};
		\node [style=root3] (91) at (2, 2) {};
		\node [style=root3] (92) at (4, 2) {};
		\node [style=root3] (93) at (6, 2) {};
		\node [style=root3] (94) at (2, 3) {};
		\node [style=root3] (95) at (4, 3) {};
		\node [style=root3] (96) at (6, 3) {};
		\node [style=root3] (97) at (2, 4) {};
		\node [style=root3] (98) at (4, 4) {};
		\node [style=root3] (99) at (6, 4) {};
	\end{pgfonlayer}
	\begin{pgfonlayer}{edgelayer}
		\draw [style=plainedge] (3.center) to (2.center);
		\draw [style=plainedge] (0.center) to (1.center);
		\foreach \i in {1,2,3,4} \draw (-0.1,\i)--(0.1,\i);
		\foreach \i in {1,2,...,7} \draw (\i,-0.1)--(\i,0.1);
	\end{pgfonlayer}
\end{tikzpicture}
}
\caption{The positive roots of the root systems of $ \mathfrak{m}_g $ when $ g $ is in the conjugacy class $ 2B $. The dark filled circles correspond to simple roots indexed by $ I_+ $. The medium filled circles on the horizontal axis correspond to the simple roots indexed by $ I_- $, i.e. the positive part of $ \mathcal{H} $. The light filled circles correspond to roots that are not simple roots. The $ \mathcal{H} $-modules can been seen in the rows above the horizontal axis with lowest-weight vectors at $ \projq{\alpha}_{(1,k)}, \projq{\alpha}_{(2,k)}, $ etc.} \label{fig:2B}
\end{figure}
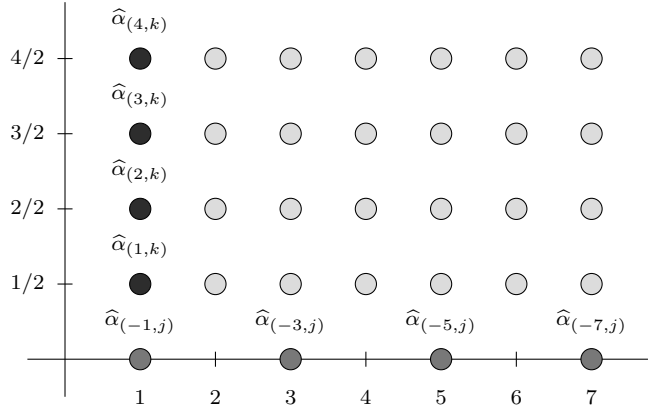
\begin{figure}[t]
\centering
{\footnotesize
\begin{tikzpicture}
	\begin{pgfonlayer}{nodelayer}
		\node [style=none] (0) at (0, 4.75) {};
		\node [style=none] (1) at (0, -0.5) {};
		\node [style=none] (2) at (7.75, 0) {};
		\node [style=none] (3) at (-0.5, 0) {};
		\node [style=root1] (4) at (2, 0) {};
		\node [style=root1] (6) at (6, 0) {};
		\node [style=root2] (12) at (1, 1) {};
		\node [style=root2] (14) at (1, 3) {};
		\node [style=root3] (21) at (3, 1) {};
		\node [style=root3] (22) at (5, 1) {};
		\node [style=root3] (27) at (3, 3) {};
		\node [style=root3] (28) at (5, 3) {};
		\node [style=none] (52) at (-0.5, 1) {$1/8$};
		\node [style=none] (53) at (-0.5, 2) {$2/8$};
		\node [style=none] (54) at (-0.5, 3) {$3/8$};
		\node [style=none] (55) at (-0.5, 4) {$4/8$};
		\node [style=none] (60) at (1, -0.5) {$1$};
		\node [style=none] (61) at (2, -0.5) {$2$};
		\node [style=none] (62) at (3, -0.5) {$3$};
		\node [style=none] (63) at (4, -0.5) {$4$};
		\node [style=none] (68) at (2, 0.5) {$\projq{\alpha}_{(-2,j)}$};
		\node [style=none] (70) at (6, 0.5) {$\projq{\alpha}_{(-6,j)}$};
		\node [style=none] (72) at (1, 1.5) {$\projq{\alpha}_{(1,k)}$};
		\node [style=none] (74) at (1, 3.5) {$\projq{\alpha}_{(3,k)}$};
		\node [style=none] (76) at (5, -0.5) {$5$};
		\node [style=none] (77) at (6, -0.5) {$6$};
		\node [style=none] (78) at (7, -0.5) {$7$};
		\node [style=root3] (80) at (7, 1) {};
		\node [style=root3] (82) at (7, 3) {};
		\node [style=root3] (85) at (4, 2) {};
		\node [style=root3] (86) at (6, 2) {};
		\node [style=root3] (89) at (2, 4) {};
		\node [style=root3] (90) at (4, 4) {};
		\node [style=root3] (91) at (6, 4) {};
		\node [style=root3] (92) at (2, 2) {};
	\end{pgfonlayer}
	\begin{pgfonlayer}{edgelayer}
		\draw [style=plainedge] (3.center) to (2.center);
		\draw [style=plainedge] (0.center) to (1.center);
		\foreach \i in {1,2,3,4} \draw (-0.1,\i)--(0.1,\i);
		\foreach \i in {1,2,...,7} \draw (\i,-0.1)--(\i,0.1);
	\end{pgfonlayer}
\end{tikzpicture}
}
\caption{The positive roots of the root systems of $ \mathfrak{m}_g $ when $ g $ is in the conjugacy class $ 4D $. The simple and non-simple roots have the same coloring system as in Figure \ref{fig:2B}. Note that there is a root at $ (2, 2/8) $ since the root space at $ (1, 1/8) $ has dimension greater than $ 1 $ and is contained in the space $ \mathcal{V}^+ $ that freely generates $ \mathfrak{u}^+ $.} \label{fig:4D}
\end{figure}
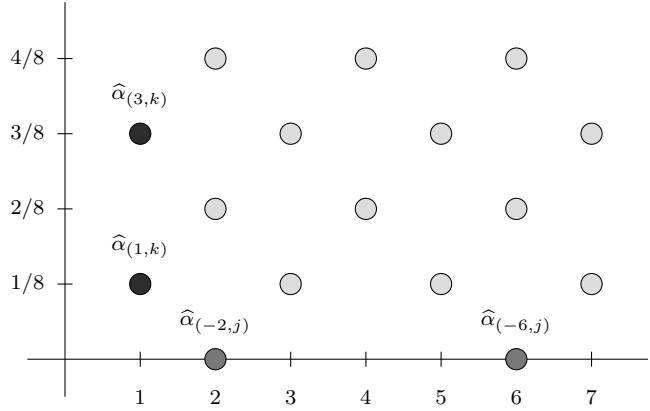

\section{Twisted denominator identities} \label{section:twisted denominator}

We can use the subalgebra $ \mathfrak{u}^+ $ to compute the denominator formula for $ \mathfrak{m}_g $. Since $ \mathfrak{g}_{I_-}^+ $ is abelian, $ \mathfrak{u}^+ $ contains enough information of $ \mathfrak{g}(A_g)^+ \cong \mathfrak{g}_{I_-}^+ \oplus \mathfrak{u}^+ $ to compute the denominator formula. Furthermore, $ \mathfrak{u}^+ $ is a free Lie algebra, making the computation of its homology immediate. The ideas in this section were introduced in \cite{JLW} for the Monster Lie algebra. These ideas were able to be followed directly in \cite{CarFricke} for the Monstrous Lie algebras of Fricke type using Theorem 5.1 of \cite{JurJPAA} and $ \mathfrak{sl}_2 $-representation theory. We follow a similar stream of ideas for the non-Fricke Monstrous Lie algebras, but with Theorem \ref{theorem:decomposition} and $ \heis $-representation theory instead. 

In \cite{CarThesis}, Carnahan computes the twisted denominator identities by computing the homology of $ \mathfrak{g}(A_g)^+ $ in the triangular decomposition of $ \mathfrak{m}_g $. Following the method in \cite{BoInvent}, he requires an extension of the results from \cite{GarLep} for homology of Kac-Moody Lie algebras to the case of Borcherds algebras. Our use of the free subalgebra $ \mathfrak{u}^+ $ avoids the use of such theory. The simplicity in this approach comes from the homology of $ \mathfrak{u}^+ $ expressed in terms of the $ \mathbb{Z} \times \frac{1}{N}\mathbb{Z} $-graded Lie algebra $ \mathcal{H}$-module $ \mathcal{V}^+ $. 

In the following, we will omit some definitions and proofs that can be found in Section 5 of \cite{JLW}.

The main observation is to apply the Euler-Poincar\'e identity
\begin{equation}\label{eq:Euler-Poincare}
H(\mathfrak{g}) = \Lambda_{-1}(\mathfrak{g}) 
\end{equation}
to $ \mathfrak{u}^+ $. Since $ \mathfrak{u}^+ $ is freely-generated by $ \V^+ $, as defined in \eqref{eq:V+-}, its homology is straight-forwardly computed to be
\begin{equation}\label{eq:homology of u}
H(\mathfrak{u}^+) = 1 - \V^+.
\end{equation}
Since $ \V^+ $ is a direct sum of $ \mathbb{Z} \times \frac{1}{N}\mathbb{Z} $-graded lowest-weight $ \heis $-modules we can compute its graded structure just from the simple roots to be
\begin{equation} \label{eq:V iso}
\V^+_{1+k, n/N} \cong \coprod_{(1^{m_1}2^{m_2} \cdots ) \vdash k} \bigotimes_{i=1}^k \Sym^{m_i} (\mathfrak{g}_{i,0} )\otimes \mathfrak{g}_{1,n/N} 
\end{equation}
Here, $ (1^{m_1}2^{m_2} \cdots ) \vdash k $ denotes a partition of $ k $ and $ \Sym^{m} $ denotes a symmetric power. This isomorphism is a linear categorification of the generating function 
\begin{align*}
\sum_{k \geq 0} \dim \V^+_{1+k, n/N}  p^k  &= \prod_{ m > 0} ( 1- p^m)^{- \dim  \mathfrak{g}_{m,0} } \dim \mathfrak{g}_{1,n/N} \\
&= \sum_{k \geq 0} \sum_{(1^{m_1}2^{m_2} \cdots ) \vdash k} \prod_{i=1}^k \binom{\dim \mathfrak{g}_{i,0} + m_i - 1}{m_i}\dim  \mathfrak{g}_{1,n/N} p^k
\end{align*}
for the dimensions of $ \V^+_{1+k, n/N} $. Furthermore, if $ \mathfrak{m}_g $ has a group $ G $ of  $ \mathbb{Z} \times \frac{1}{N}\mathbb{Z} $-graded Lie algebra automorphisms, then \eqref{eq:V iso} still holds as an isomorphism of representations of $ G $. For example, this group can be taken to be a certain central extension of the centralizer of $ g $ in $ \M $ as shown in \cite{CarIV}. 

Denoting the grading with the powers of $ p $ and $ q $ gives us
\begin{equation}
H(\mathfrak{u}^+) = 1 - \sum_{\substack{ k \geq 0 \\ n > 0}} \coprod_{(1^{m_1}2^{m_2} \cdots ) \vdash k} \bigotimes_{i=1}^k \Sym^{m_i} (\mathfrak{g}_{i,0} )\otimes \mathfrak{g}_{1,n/N} p^{1+k}q^{n/N}.
\end{equation}
We also have
\begin{equation}
\Lambda_{-1}(\mathfrak{u}^+) = \exp \left( -\sum_{k=1}^\infty \frac{1}{k} \Psi^k (\mathfrak{u}^+) \right) = \exp \left( -\sum_{k=1}^\infty \frac{1}{k} \sum_{m,n > 0}  \Psi^k(\mathfrak{g}_{m,n/N} ) p^{km} q^{kn/N} \right) ,
\end{equation}
where $ \Psi^k $ are Adams operations (see \cite{JLW} for the details). We can take the trace of $ h \in G $ acting on \eqref{eq:Euler-Poincare} to obtain twisted denominator identities for $ \mathfrak{m}_g $. We have
\begin{align*}
&\Tr(h \vert H(\mathfrak{u}^+)) \\
&\quad = 1 - \sum_{\substack{ k \geq 0 \\ n > 0}} \sum_{(1^{m_1}2^{m_2} \cdots ) \vdash k} \prod_{i=1}^k \Tr(h \vert \Sym^{m_i} (\mathfrak{g}_{i,0} )) \Tr(h \vert \mathfrak{g}_{1,n/N}) p^{1+k}q^{n/N}\\
&\quad =  1 - \sum_{k \geq 0} \sum_{(1^{m_1}2^{m_2} \cdots ) \vdash k} \prod_{i=1}^k \Tr(h \vert \Sym^{m_i} (\mathfrak{g}_{i,0} )) p^k  \sum_{n > 0} \Tr(h \vert \mathfrak{g}_{1,n/N}) p q^{n/N} \\
&\quad =  1 - \prod_{m > 0} \sum_{j \geq 0} \Tr(h \vert \Sym^{m}(\mathfrak{g}_{i,0} ) )(p^m)^j  \sum_{n > 0} \Tr(h \vert \mathfrak{g}_{1,n/N}) p q^{n/N} \\
&\quad =  1 - \prod_{m > 0} \exp \left( \sum_{k>0} \frac{1}{k} \Tr( h^k \vert \mathfrak{g}_{m,0} ) (p^{m})^k \right)  \sum_{n > 0} \Tr(h \vert \mathfrak{g}_{1,n/N}) p q^{n/N} ,
\end{align*}
where the last step uses the generating function for characters of symmetric powers of representations of $ G $. We also have
\begin{align*}
\Tr(h \vert \Lambda_{-1}(\mathfrak{u}^+)) &= \exp \left( -\sum_{k>0} \frac{1}{k} \sum_{m,n > 0}  \Tr( h^k \vert \mathfrak{g}_{m,n/N} ) p^{km} q^{kn/N} \right) \\
&= \prod_{m,n > 0} \exp \left( -\sum_{k>0} \frac{1}{k} \Tr( h^k \vert \mathfrak{g}_{m,n/N} ) (p^{m} q^{n/N})^k \right) \\
&= \prod_{\substack{m > 0\\ n\geq 0}} \exp \left( -\sum_{k>0} \frac{1}{k} \Tr( h^k \vert \mathfrak{g}_{m,n/N} ) (p^{m} q^{n/N})^k \right)  \\
&\qquad \cdot \prod_{m > 0} \exp \left( \sum_{k>0} \frac{1}{k} \Tr( h^k \vert \mathfrak{g}_{m,0} ) (p^{m})^k \right).
\end{align*}
These two expressions can be equated to give the twisted denominator identity
\begin{multline}
\prod_{m > 0} \exp \left( - \sum_{k>0} \frac{1}{k} \Tr( h^k \vert \mathfrak{g}_{m,0} )  (p^{m})^k \right) - \sum_{n > 0} \Tr(h \vert \mathfrak{g}_{1,n/N}) p q^{n/N} \\
= \prod_{\substack{m > 0\\ n\geq 0}} \exp \left( -\sum_{k>0} \frac{1}{k} \Tr( h^k \vert \mathfrak{g}_{m,n/N} ) (p^{m} q^{n/N})^k \right) .
\end{multline}
The twisted denominator identities computed here match those as in \cite{CarThesis}. 

When $ h = \id $, we obtain the denominator formula
\begin{equation}
 \prod_{m > 0}(1-p^m)^{c(m,0)} -  \sum_{n >0} c(1,n/N) p q^{n/N} =  \prod_{\substack{m > 0 \\ n \geq 0}} (1-p^m q^{n/N})^{c(m,n/N)}.
\end{equation}
Multiplying this indentity by $ p^{-1} $ and recalling the definition of $ \{c(m,0)\}_{m =1}^\infty $ and  $ \{c(1,n/N)\}_{n =1}^\infty $ from \eqref{eq:c(m,0) defn} and \eqref{eq:c(1,n/N) defn}, we get
\begin{equation}\label{eq:denomIDnon-fricke}
T_g(\sigma) - T_g(-1/\tau) = p^{-1} \prod_{\substack{m > 0 \\ n \geq 0}} (1-p^m q^{n/N})^{c(m,n/N)},
\end{equation}
where $ p = e^{2 \pi i \sigma} $ and $ q = e^{2 \pi i \tau} $. This recovers the generating functions used to define the Monstrous Lie algebras.

\bibliographystyle{amsalpha}
\bibliography{monstref.bib} 
	
\end{document}